\RequirePackage[leqno]{amsmath}
\documentclass[leqno,a4paper]{amsart}
\usepackage[usenames,dvipsnames]{xcolor}
\usepackage{amsfonts}
\usepackage{amsmath}
\usepackage{enumitem}
\setlist[itemize]{topsep=0pt,after=\vspace{1.5\baselineskip}}
 \usepackage[colorlinks=true]{hyperref}
\usepackage{empheq}
\usepackage[T1]{fontenc}
\usepackage[utf8]{inputenc}
\usepackage{todonotes}

\def\R{\mathbb R}  
\newtheorem{theorem}{Theorem}[section]
\newtheorem{corollary}[theorem]{Corollary}
\newtheorem{lemma}[theorem]{Lemma}

\newcommand{\na}{\nabla}
\newcommand{\Om}{\Omega}
\newcommand{\Ombar}{\overline{\Omega}}

\newcommand{\f}[2]{\frac{#1}{#2}}
\newcommand{\io}{\int_\Omega}

\newcommand{\norm}[2][]{\left\|#2\right\|_{#1}}

\newcommand{\Tmax}{T_{max}}

\usepackage{newunicodechar}
\newunicodechar{∞}{\infty}
\newunicodechar{α}{\alpha}
\newunicodechar{β}{\beta}
\newunicodechar{γ}{\gamma}
\newunicodechar{δ}{\delta}
\newunicodechar{ε}{\varepsilon}
\newunicodechar{κ}{\kappa}
\newunicodechar{χ}{\chi}
\newunicodechar{μ}{\mu}
\newunicodechar{σ}{\sigma}
\newunicodechar{τ}{\tau}
\newunicodechar{ξ}{\xi}
\newunicodechar{φ}{\varphi}
\newunicodechar{Δ}{\Delta}
\newunicodechar{ℝ}{\mathbb{R}}

\title[Analysis of a singular Keller-Segel model with consumption] 
      {Finite-time blow-up prevention in a two-dimensional chemotaxis-consumption model with weakly singular sensitivity}
\author[G. Viglialoro]{}
\keywords{Nonlinear parabolic systems, chemotaxis,  singular sensitivity, global existence.\\
\textit{Acknowledgments:} GV is member of the Gruppo Nazionale per l'Analisi Matematica, la Probabilit\`a e le loro Applicazioni (GNAMPA) of the Istituto Na\-zio\-na\-le di Alta Matematica (INdAM) and is partially supported by the research project \textit{Integro-differential Equations and Non-Local Problems}, funded by Fondazione di Sardegna (2017). }
\RequirePackage[normalem]{ulem} 
\RequirePackage{color}\definecolor{RED}{rgb}{1,0,0}\definecolor{BLUE}{rgb}{0,0,1} 
\providecommand{\DIFaddbegin}{} 
\providecommand{\DIFaddend}{} 
\providecommand{\DIFdelbegin}{} 
\providecommand{\DIFdelend}{} 
\providecommand{\DIFaddbeginFL}{} 
\providecommand{\DIFaddendFL}{} 
\providecommand{\DIFdelbeginFL}{} 
\providecommand{\DIFdelendFL}{} 
\newcommand{\DIFscaledelfig}{0.5}
\RequirePackage{settobox} 
\RequirePackage{letltxmacro} 
\newsavebox{\DIFdelgraphicsbox} 
\newlength{\DIFdelgraphicswidth} 
\newlength{\DIFdelgraphicsheight} 
\LetLtxMacro{\DIFOincludegraphics}{\includegraphics} 
\newcommand{\DIFaddincludegraphics}[2][]{{\color{blue}\fbox{\DIFOincludegraphics[#1]{#2}}}} 
\newcommand{\DIFdelincludegraphics}[2][]{
\sbox{\DIFdelgraphicsbox}{\DIFOincludegraphics[#1]{#2}}
\settoboxwidth{\DIFdelgraphicswidth}{\DIFdelgraphicsbox} 
\settoboxtotalheight{\DIFdelgraphicsheight}{\DIFdelgraphicsbox} 
\scalebox{\DIFscaledelfig}{
\parbox[b]{\DIFdelgraphicswidth}{\usebox{\DIFdelgraphicsbox}\\[-\baselineskip] \rule{\DIFdelgraphicswidth}{0em}}\llap{\resizebox{\DIFdelgraphicswidth}{\DIFdelgraphicsheight}{
\setlength{\unitlength}{\DIFdelgraphicswidth}
\begin{picture}(1,1)
\thicklines\linethickness{2pt} 
{\color[rgb]{1,0,0}\put(0,0){\framebox(1,1){}}}
{\color[rgb]{1,0,0}\put(0,0){\line( 1,1){1}}}
{\color[rgb]{1,0,0}\put(0,1){\line(1,-1){1}}}
\end{picture}
}\hspace*{3pt}}} 
} 
\LetLtxMacro{\DIFOaddbegin}{\DIFaddbegin} 
\LetLtxMacro{\DIFOaddend}{\DIFaddend} 
\LetLtxMacro{\DIFOdelbegin}{\DIFdelbegin} 
\LetLtxMacro{\DIFOdelend}{\DIFdelend} 
\DeclareRobustCommand{\DIFaddbegin}{\DIFOaddbegin \let\includegraphics\DIFaddincludegraphics} 
\DeclareRobustCommand{\DIFaddend}{\DIFOaddend \let\includegraphics\DIFOincludegraphics} 
\DeclareRobustCommand{\DIFdelbegin}{\DIFOdelbegin \let\includegraphics\DIFdelincludegraphics} 
\DeclareRobustCommand{\DIFdelend}{\DIFOaddend \let\includegraphics\DIFOincludegraphics} 
\LetLtxMacro{\DIFOaddbeginFL}{\DIFaddbeginFL} 
\LetLtxMacro{\DIFOaddendFL}{\DIFaddendFL} 
\LetLtxMacro{\DIFOdelbeginFL}{\DIFdelbeginFL} 
\LetLtxMacro{\DIFOdelendFL}{\DIFdelendFL} 
\DeclareRobustCommand{\DIFaddbeginFL}{\DIFOaddbeginFL \let\includegraphics\DIFaddincludegraphics} 
\DeclareRobustCommand{\DIFaddendFL}{\DIFOaddendFL \let\includegraphics\DIFOincludegraphics} 
\DeclareRobustCommand{\DIFdelbeginFL}{\DIFOdelbeginFL \let\includegraphics\DIFdelincludegraphics} 
\DeclareRobustCommand{\DIFdelendFL}{\DIFOaddendFL \let\includegraphics\DIFOincludegraphics} 
\begin{document}
\maketitle
\DIFdelbegin 
\DIFdelend 

\centerline{\scshape Giuseppe Viglialoro}
\medskip
{
 \footnotesize
 \centerline{Dipartimento di Matematica e Informatica}
 \centerline{Universit\`{a} di Cagliari}
 \centerline{V. le Merello 92, 09123. Cagliari (Italy)}
  \centerline{giuseppe.viglialoro@unica.it}
 \medskip
}

\bigskip
\begin{abstract}
In this paper we study a two-dimensional chemotaxis-consumption system 
with singular sensitivity  and endowed with Neumann boundary conditions. Sufficient conditions on the data of the problem are given so that the globability of classical solutions is shown, thus excluding any finite-time blow-up scenario. 
%
%
\end{abstract}
\section{Introduction, motivation and claim of the main theorem}\label{IntroductionSection} 
This paper deals with a variant of the classical Keller--Segel system (\cite{K-S-1970}) which models chemotaxis  phenomena of cells at the position $x$ of their environment  and at the time $t$ when their motion is influenced by a chemical signal. We are interested in the situation where the distribution $u=u(x,t)$ of such cells, occupying an insulated domain, directs their movement in response to a substance $v=v(x,t)$ they consume, once  the initial distribution of the cells and of the chemical concentration are known.

In particular, the mathematical formulation of the biological model described above is the initial-boundary value problem
\begin{equation}\label{problem}
\begin{cases}
u_t=\Delta u -\chi \nabla \cdot \big(\frac{u}{v^\gamma} \nabla v\big)   & \text{in } \Omega\times(0,\infty), \\
v_t=\Delta v-f(u) v & \text{in } \Omega\times(0,\infty),\\
\frac{\partial u}{\partial \nu}=\frac{\partial v}{\partial \nu}=0 & \text{in } \partial\Omega\times(0,\infty),\\
u(x,0)=u_{0}(x) \quad \textrm{and}\quad v(x,0)=v_0(x),& x\in  \Omega,
\end{cases}
\end{equation}
where $\chi>0$,
 \begin{equation}\label{assumptions}\tag{$\mathcal{A}$} 
 \begin{cases}
  \Omega\subset ℝ^2 \text{ is a smooth and bounded domain},\\
\beta\in (0,1] \quad\text{and}\quad 0<\gamma <1, \\
  f\in C^1(ℝ) \quad\text{and}\quad 0\le f(s)\le s^\beta \quad \text{for all } s>0,
 \end{cases}
\end{equation}
 and 
\begin{equation}\label{initdata}
(u_0,v_0)\in C^0(\bar{\Omega})\times W^{1,r}(\Omega) \text{ for some } r>2, \text{ satisfy } u_0\ge 0 \text{ and } v_0>0 \text{ in  }\Ombar.
\end{equation}
The aim of this investigation is to extend \cite[Theorem 2.1]{LankeitViglialoro}, where in \eqref{problem} it is $\gamma=1$, so that the sensitivity is reduced to $u/v$: \textit{For any $\chi\in (0,1)$ and $\beta\not\equiv 1$ in \eqref{assumptions}, and any initial data as in \eqref{initdata}, the existence of global classical solutions is established.} In accordance to this result, our main theorem shows that weakening the singularity of the sensitivity $u/v$ at $v=0$ to $u/v^\gamma$ suffices to guarantee the same conclusions of \cite[Theorem 2.1]{LankeitViglialoro} even for arbitrarily large values of $\chi$ provided that smallness assumptions on the maximum of $v_0$ on $\bar{\Omega}$ are satisfied.

To be precise we will prove the following
\begin{theorem}\label{MainTheorem}   
Let  \eqref{assumptions} with $\beta \not\equiv 1$. Then for any $\chi>0$ and $(u_0,v_0)$ as in \eqref{initdata} satisfying 
$\lVert v_0 \rVert_{L^\infty(\Omega)} < \chi^\frac{1}{\gamma-1}$, there exists a unique pair of functions $(u,v)$, 
\begin{align}\label{solutionregularity}
\begin{cases}
 u\in C^0(\bar{\Omega}\times [0,\infty))\cap C^{2,1}(\bar{\Omega}\times (0,\infty)),\\
 v\in C^0(\bar{\Omega}\times [0,\infty))\cap C^{2,1}(\bar{\Omega}\times (0,\infty))\cap L^\infty_{loc}([0,\infty), W^{1,r}(\Om)), \end{cases}
\end{align}
which solve problem \eqref{problem}.
\end{theorem}
\subsection{Some preparatory inequalities} We will rely on these two lemmas. 
 \begin{lemma}[Gagliardo-Nirenberg and Neumann heat semigroup inequalities]\label{InequalityG-NLemma}
 Let $\Om\subset ℝ^2$ be a smooth and bounded domain, $0< \beta\leq 1,$ $ p\in [2,4]$ and $\theta=1-\frac{2}{p}\in [0,\frac{1}{2}].$
 Then there exist $C_{GN}=C_{GN}(\Omega)>0$ and $C_S=C_S(\Omega,\beta)>0$ such that 
\item  []  
 \begin{equation}\label{InequalityTipoG-N} 
 \|  \varphi \|_{L^{p}(\Omega)} \leq C_{GN} ( \| \nabla  \varphi \|_{L^2(\Omega)}^{\theta} \|  \varphi \|_{L^2(\Omega)}^{1 - \theta}+  \|  \varphi \|_{L^2(\Omega)}) \quad \forall \;  \varphi\in H^1(\Omega);
 \end{equation}
 \item [] 
 \begin{equation}\label{LpLqEstimateGradient0}   
 \lVert  e^{t\Delta}  \varphi  \lVert_{L^\infty(\Omega)}\leq C_{S} (1+t^{-\frac{\beta}{2}}) \lVert  \varphi\lVert_{L^\frac{2}{\beta}(\Omega)}\quad \forall \;  \varphi\in  L^\frac{2}{\beta}(\Omega).
 \end{equation}
 \begin{proof}
 See \cite[p. 126]{Nirenber_GagNir_Ineque} and  \cite[Lemma 1.3]{WinklAggre},  respectively.
 \end{proof}
 \end{lemma}
\begin{lemma}\label{lem:estimate:L2:spacetime}
 Let $\Om\subset ℝ^2$ be a smooth and bounded domain, let $T>0$, $L_1>0$, $m>0$ and $C_{GN}$ the constant introduced in Lemma \ref{InequalityG-NLemma}. If $ 0\leq \varphi\in C^0(\Ombar\times[0,T))\cap C^{2,1}(\Ombar\times(0,T))$ fulfills 
\begin{equation}\label{EstimatesMassEnGradientOverFunction}
 \io  \varphi(\cdot,t)=m \quad \text{and } \quad \int_0^t\io \f{|\na  \varphi|^2}{ \varphi } \le L_1(1+t) \qquad \text{for all } t\in(0,T),
\end{equation}
then for  $\mathcal{K}=2mC_{GN}^4(L_1+4m)$ we have
\begin{equation}\label{GagliardoNirenberg_u^2_on(0,T)}
 \int_0^t \io  \varphi^2 \le \mathcal{K}(1+t) \quad \textrm{ for all } t<T.
\end{equation}
\begin{proof}
Inequality \eqref{InequalityTipoG-N} and $(A+B)^4\leq 8(A^4+B^4)$, $A,B\in \R$, enable us to estimate
 \begin{equation*} 
 \begin{split}
 \int_\Omega  \varphi^{2}=\lVert \sqrt{ \varphi}\rVert^{4}_{L^{4}(\Omega)} &\leq 
 8C_{GN}^{4}\Big[\frac{m}{4}\Big(\int_\Omega \frac{\lvert \nabla  \varphi \rvert^2}{ \varphi}\Big)+m^{2}\Big] \qquad \text{ on } (0,T),
 \end{split}
 \end{equation*} 
 and we conclude by virtue of the assumption on $\int_0^t\io \f{|\na  \varphi|^2} \varphi$. 
\end{proof}
\end{lemma}
\section{Proof of Theorem \ref{MainTheorem} and one consequence}\label{sec:globalexistence}
\subsection{Toward the proof of the main theorem} Let us first claim a result on local-in-time existence of classical solutions to system \eqref{problem}, whose proof  is standard  (see, e.g., \cite[Lemma 3.1]{BellomoEtAl} or \cite[Lemma 2.2]{win_ct_sing_abs}).   
\begin{lemma}\label{LocalExistenceLemma}
Assume \eqref{assumptions}. 
For  any $\chi>0$ and $(u_0,v_0)$ as in \eqref{initdata}, there are $T_{max}\in (0,∞]$ and a uniquely determined pair of functions $(u,v)$ with regularity as in \eqref{solutionregularity} which solve problem \eqref{problem} in $\Om\times(0,T_{max})$ and are such that if $T_{max}<\infty$ then
\begin{equation}\label{extensibility_criterion_Eq} 
\limsup_{t\nearrow T_{max}}(\lVert u (\cdot,t)\rVert_{L^\infty(\Omega)}+\lVert v (\cdot,t)\rVert_{W^{1,r}(\Omega)})=\infty.
\end{equation}
Moreover, we have $u\geq 0$ and $0<v\leq \lVert v_0 \rVert_{L^\infty(\Omega)}$ in $\bar{\Omega} \times (0,T_{max})$, as well as
\begin{equation}\label{Bound_of_u} 
\int_\Omega u (\cdot ,t)  =m=\int_\Omega u_0 \quad \textrm{for all}\quad  t\in (0,T_{max}).
\end{equation}
%
%
\end{lemma}
With the local solution to \eqref{problem} in our hands, the suitable transformation (see e.g. \cite{LankeitLocallyBoundedSingularity,win_ct_sing_abs})
\begin{equation}\label{transformationVtoW}
w:=-\log \Big(\frac{v}{\lVert v_0 \rVert_{L^\infty(\Omega)}}\Big), \qquad w_0:=-\log \Big(\frac{v_0}{\lVert v_0 \rVert_{L^\infty(\Omega)}}\Big), 
\end{equation}
implies $w\geq 0$ in $\Omega \times (0,T_{max})$ and that $(u,w)$ also classically solves this problem
\begin{equation}\label{AuxiliaryMainsSystem}
\begin{cases}
u_t=\Delta u +\chi \lVert v_0\rVert_{L^\infty(\Omega)}^{1-\gamma}\nabla \cdot (ue^{-(1-\gamma)w} \nabla w) &  \text{ in }\Omega \times (0,T_{max}), \\
w_t=\Delta w -\rvert \nabla w \rvert^2+f(u) &  \text{ in }\Omega \times (0,T_{max}), \\
\frac{\partial u}{\partial \nu}=\frac{\partial w}{\partial \nu}=0 & \text{ in }\partial \Omega\times (0,T_{max}),\\
u(x,0)=u_{0}(x)\geq 0 \quad w(x,0)=w_0(x)\geq 0,& x\in  \bar{\Omega}.
\end{cases}
\end{equation}
Such transformed problem allows us, inter alia, to sharpen the extensibility criterion \eqref{extensibility_criterion_Eq} so to warrant the globability of local solutions. 
\begin{lemma}\label{LemmaToEnsureTisInfty} 
Assume \eqref{assumptions}. For any $\chi>0$ and $(u_0,v_0)$ as in \eqref{initdata}, let $(u,v)$ be the local-in-time classical solution of problem \eqref{problem} provided by Lemma \ref{LocalExistenceLemma}, and $(u,w)$ that of problem \eqref{AuxiliaryMainsSystem}, with $w$ as in \eqref{transformationVtoW}. If there exists $C>0$ such that 
\begin{equation}\label{Bounds_of_SqrtU_And_w} 
\begin{cases}
\int_0^t\int_\Omega\f{|\na u|^2}{u } \leq C(1+t)&\forall \;t\in(0,T_{max}),\\
 \int_\Omega u(\cdot, t)\log u(\cdot,t)\leq C(1+t)& \forall \;t\in(0,T_{max}),\\
 \lVert w (\cdot,t)\rVert_{L^\infty(\Omega)}\leq C(1+t)& \forall \;t\in(0,T_{max}),
\end{cases}
\end{equation}
then $T_{max}=\infty.$
\begin{proof}
If $T_{max}$ were finite, from \eqref{transformationVtoW} we would have that 
$1/v=e^w/\lVert v_0\rVert_{L^\infty(\Omega)},$
and the third assumption in \eqref{Bounds_of_SqrtU_And_w} would guarantee that 
 $1/v\le L$ in $\Om\times(0,\Tmax),$ with some $L>0$.  For such a constant $L$, we introduce a smooth \textit{cutoff function} $\xi_L:\R \rightarrow [0,1]$, decreasing and verifying $\xi_L(v)=1$ for $v\leq 1/(2L)$ and $\xi_L(v)=0$ for $v\geq 1/L$. Subsequently, 
for any $(x,t,u,v) \in \bar{\Omega}\times [0,\infty)\times \R^2$ and for some $\omega \in (0,1)$, 
 \begin{equation*}
 \begin{split}
 S_\gamma(x,t,u,v)=
 \xi_L(v)\frac{2\chi}{L}+(1-\xi_L(v))\frac{\chi}{v^\gamma} \in C^{1+\omega}_{\textrm{loc}}(\bar{\Omega}\times [0,\infty)\times \mathbb{R}^2), 
 \end{split}
 \end{equation*}
and additionally satisfies $S_\gamma(x,t,u,v)\equiv \chi/v^\gamma$ for all $v\geq 1/L.$ Hence, coherently with the nomenclature of \cite[Sec. 3]{BellomoEtAl}, setting
\begin{align*}
f(x,t,u,v)&\equiv 0,\quad \textrm{and}\quad g(x,t,u,v)=v-f(u) v,
\end{align*}
the two partial differential equations of problem \eqref{problem} read 
 \[
 \begin{cases}
  u_t=\Delta u-\nabla \cdot (uS_\gamma(x,t,u,v)\nabla v)+f(x,t,u,v)& \text{in }\Omega \times (0,T_{max})\\
  v_t=\Delta v-v+g(x,t,u,v) & \text{in }\Omega \times (0,T_{max}).
 \end{cases}
\]
In particular, we have $S_\gamma\in C^{1+\omega}_{\textrm{loc}}(\bar{\Omega}\times [0,\infty)\times \mathbb{R}^2)$,  for any $0<\gamma<1$, $f\in\ C^{1-}_{\textrm{loc}}(\bar{\Omega}\times [0,\infty)\times \mathbb{R}^2)$ and $g\in C^{1-}_{\textrm{loc}}(\bar{\Omega}\times [0,\infty)\times \mathbb{R}^2)$,  $f(x,t,0,v)=0$ for all $(x,t,v)\in  \bar{\Omega}\times [0,\infty)\times \mathbb{R}$ and $g(x,t,u,0 )=0$ for all $(x,t,u)\in \bar{\Omega}\times [0,\infty)\times \mathbb{R}$. 

On the other hand, the Young inequality, the second equation of \eqref{problem}, the estimate on $f$ given in \eqref{assumptions} and the upper bound for $v$ in Lemma \ref{LocalExistenceLemma} provide
 \begin{equation*}
 \begin{split}
 \frac{d}{dt}\int_\Omega \lvert \nabla v\rvert^2 &=2\int_\Omega \nabla v \cdot \nabla v_t\leq 
  -2\int_\Omega (\Delta v)^2+\int_\Omega  (\Delta v)^2+\int_\Omega  u^{2\beta}v^2
 \\&
 \leq  -\int_\Omega (\Delta v)^2+ \lVert v_0 \rVert_{L^\infty(\Omega)}^2\beta \int_\Omega u^2 + \lVert v_0 \rVert_{L^\infty(\Omega)}^2(1-\beta)\lvert \Omega\rvert \quad \text{in } (0,\Tmax),
 \end{split}
 \end{equation*}
and, neglecting the nonpositive term  $-\int_\Omega (\Delta v)^2$, we obtain 
 \begin{equation*}
 \frac{d}{dt}\int_\Omega \lvert \nabla v\rvert^2 \leq \lVert v_0 \rVert_{L^\infty(\Omega)}^2\beta \int_\Omega u^2 +\lVert v_0 \rVert_{L^\infty(\Omega)}^2(1-\beta)\lvert \Omega\rvert \qquad \text{in } (0,\Tmax).
 \end{equation*}
 Finally, by means of  \eqref{GagliardoNirenberg_u^2_on(0,T)} with $\varphi=u$, an integration over $(0,t)$ in the previous step produces
 \begin{equation}\label{BoundSquareNablav(0,T)}
 \int_\Omega \lvert \nabla v(\cdot,t)\rvert^2 \leq C_2(1+t)\quad \textrm{for all}\quad t\in (0,T_{max}),
 \end{equation}
 where 
  \[
C_2=\lVert v_0 \rVert_{L^\infty(\Omega)}^2\beta \mathcal{K}+\max\bigg\{\int_\Omega |\nabla v_0|^2,\lVert v_0 \rVert_{L^\infty(\Omega)}^2(1-\beta)\lvert \Omega\rvert\bigg\}.
\]
Subsequently, in view of the assumptions in \eqref{Bounds_of_SqrtU_And_w}, estimate \eqref{BoundSquareNablav(0,T)}, the regularity and boundedness of both $S_\gamma$ and $v$ and the expression of $g$ given above, there exists a positive $N$ such that for all $t\in(0,T_{max})$ 
\begin{equation*} 
\begin{cases}
\int_\Omega \lvert \nabla v(\cdot,t)\rvert^2\leq N,\quad \int_\Omega u(\cdot, t)\log u(\cdot,t)\leq N, \\ 
S_\gamma(x,t,u,v)\leq N,\quad\lvert  g(x,t,u,v)\rvert \leq N(1+u).
\end{cases}
\end{equation*}
From this, \cite[Lemma 3.3]{BellomoEtAl} implies that $t\mapsto\lVert u(\cdot, t)\rVert_{L^\infty(\Omega)}+\lVert v(\cdot, t)\rVert_{W^{1,r}(\Omega)}$ is bounded on $(0,T_{max})$, contrasting the extensibility criterion \eqref{extensibility_criterion_Eq}: so $T_{max}=\infty.$
\end{proof}
\end{lemma}
Now, we define and study the evolution in time of this  functional 
\begin{equation}\label{FuntionalUsedForGlobalExistence}
\mathfrak{F}(t)=\mathfrak{F}(u,w):=\int_\Omega u \log u+\frac{1}{2}\int_\Omega u w\quad \textrm{on}\quad (0,T_{max}),
\end{equation}
with initial value  $\mathfrak{F}(0)=\int_\Omega u _0\log u_0+\frac{1}{2}\int_\Omega u_0 w_0.$
\begin{lemma}\label{LemmaControllinguLoguAndnablauOveru}
Assume \eqref{assumptions} with $\beta\not\equiv 1$. For any $\chi>0$ and $(u_0,v_0)$ as in \eqref{initdata} satisfying
\begin{equation}\label{ConditionOnChiAndv0}
\lVert v_0 \rVert_{L^\infty(\Omega)} < \chi^\frac{1}{\gamma-1},
\end{equation} 
let $(u,v)$ be the local-in-time classical solution of problem \eqref{problem} provided by Lemma \ref{LocalExistenceLemma}, and $(u,w)$ that of problem \eqref{AuxiliaryMainsSystem}, with $w$ as in \eqref{transformationVtoW}. Then for some  $L_1, L_2>0$, these relations hold: 
 \begin{equation}\label{FirstBoundNablaSqrtuOn(0,T)} 
 \int_0^t\int_\Omega \frac{\lvert \nabla u\rvert ^2}{u}\leq L_1(1+t) \quad\textrm{for all}\quad t \in (0,T_{max});
 \end{equation}
\begin{equation}\label{BoundOfEnergyFunctional}
\int_\Omega u (\cdot,t)\log u(\cdot,t) \leq L_2(1+t) \quad\textrm{for all}\quad t \in (0,T_{max}).
\end{equation}
\begin{proof}
From \eqref{FuntionalUsedForGlobalExistence}, $s \log s \geq -\frac{1}{e}$ for all $s> 0$, and  $u\geq 0$ and $w\geq 0$, we have
\begin{equation}\label{BoundedFromBelowF}
\mathfrak{F}(u,w)\geq -\frac{\lvert \Omega \rvert}{e}\quad \textrm{ on } (0,T_{max}).
\end{equation}
Using \eqref{AuxiliaryMainsSystem}, a differentiation of  \eqref{FuntionalUsedForGlobalExistence} and the divergence theorem infer for $t<T_{max}$
\begin{equation*}
\begin{split}
\frac{d}{d t}\mathfrak{F}& \leq-\int_\Omega \frac{\lvert \nabla u \rvert^2}{u}-\chi\lVert v_0\rVert_{L^\infty(\Omega)}^{1-\gamma}\int_\Omega e^{-(1-\gamma)w}\nabla u \cdot \nabla w- \int_\Omega \nabla u \cdot \nabla w\\ & 
\quad - \frac{\chi \lVert v_0\rVert_{L^\infty(\Omega)}^{1-\gamma}}{2}\int_\Omega u  e^{-(1-\gamma)w}\lvert \nabla w\rvert^2-\frac{1}{2}\int_\Omega u |\nabla w|^2  +\frac{1}{2}\int_\Omega u^{\beta+1},
\end{split}
\end{equation*}
where we also considered \eqref{Bound_of_u} and the bound for $f$ in \eqref{assumptions}. Then, in order to absorb the  integrals involving $|\nabla w|^2$, we use the Young inequality in those containing $\nabla u \cdot \nabla w$: since $e^{-(1-\gamma)w}\leq 1$ for $0<\gamma<1$ and $w\geq 0$ in $\Omega\times (0,T_{max})$, we obtain
 \begin{equation}\label{DerivativeF_0}
 \begin{split}
 &\frac{d}{d t}\mathfrak{F} 
 \le\bigg(-\frac{1}{2}+\frac{\chi \lVert v_0\rVert_{L^\infty(\Omega)}^{1-\gamma}}{2}\bigg)\int_\Omega \frac{\lvert \nabla u \rvert^2}{u} +\frac{1}{2}\int_\Omega u^{\beta+1}\; \textrm{on}\;\quad  (0,T_{max}),
 \end{split}
 \end{equation}
so that due to the hypothesis \eqref{ConditionOnChiAndv0}, $c_0:=\frac{1}{2}-\frac{\chi \lVert v_0\rVert_{L^\infty(\Omega)}^{1-\gamma}}{2}>0$ and \eqref{DerivativeF_0} reads
 \begin{equation}\label{DerivativeF_1}
 \begin{split}
 &\frac{d}{d t}\mathfrak{F} 
 \le-c_0\int_\Omega \frac{\lvert \nabla u \rvert^2}{u} +\frac{1}{2}\int_\Omega u^{\beta+1}\; \textrm{on}\; (0,T_{max}).
 \end{split}
 \end{equation}
As to $\frac{1}{2}\int_\Omega u^{\beta+1}$,  \eqref{InequalityTipoG-N}  with $\varphi=\sqrt{u}$ and  the  mass conservation property \eqref{Bound_of_u} provide 
\begin{equation}\label{GagliardoNirenberBoundu^beta+1}
\begin{split}
\int_\Omega u^{\beta+1}&=\lVert \sqrt{u}\rVert^{2(\beta+1)}_{L^{2(\beta+1)}(\Omega)}\leq 
 2^{2\beta+1}C_{GN}^{2(\beta+1)}\Big[\frac{m}{4^\beta}\Big(\int_\Omega \frac{\lvert \nabla u \rvert^2}{u}\Big)^\beta+m^{\beta+1}\Big]
\end{split}
\end{equation} 
on the interval $(0,T_{max})$, where $(A+B)^k\leq 2^{k-1}(A^k+B^k)$, $A,B\geq 0 \textrm{ and } k\geq 1$ has been also employed.
In account of Young's inequality
\begin{equation}\label{GagliardoNirenberBoundu^beta+1Young}
\begin{split}
\int_\Omega u^{\beta+1}
\leq c_0\int_\Omega \frac{\lvert \nabla u \rvert^2}{u}+c_1m^\frac{1}{1-\beta}+ 2^{2\beta+1}C_{GN}^{2(\beta+1)} m^{\beta+1}\quad \textrm{ on } (0,\Tmax),
\end{split}
\end{equation} 
where $c_1=2(1-\beta) C_{GN}^{\frac{2(\beta+1)^2}{1-\beta}}(\frac{c_0}{2 \beta })^\frac{\beta}{\beta-1}$.  Hence, by inserting \eqref{GagliardoNirenberBoundu^beta+1Young} into \eqref{DerivativeF_1} and setting $c_2=\frac{c_1}{2}+4^\beta C_{GN}^{2(\beta+1)} m^{\beta+1}$ we get
\begin{equation}\label{InequalityForF}
\frac{d}{d t}\mathfrak{F} +\frac{c_0}{2}\int_\Omega \frac{\lvert \nabla u\rvert^2}{u}\leq c_2 \quad \textrm{on}\quad  (0,T_{max}),
\end{equation}
and an integration on $(0,t)$, for $t < T_{max}$, together with expression \eqref{BoundedFromBelowF}, gives
\begin{equation*}
-\frac{\lvert \Omega \rvert}{e}\leq \mathfrak{F}(t) =\int_\Omega u(\cdot,t)\log u(\cdot,t)\leq \mathfrak{F}(0) +c_2 t \quad \textrm{for all}\quad t\in (0,T_{max}).
\end{equation*}
In view of this bound, another integration of \eqref{InequalityForF} on $(0,t)$ yields 
\[  
\frac{c_0}{2}\int_0^t\int_\Omega \frac{\lvert \nabla u\rvert^2}{u}\leq c_2t+\mathfrak{F}(0)+\frac{\lvert \Omega \rvert}{e} \qquad \text{for } t\in(0,\Tmax),
\]
and \eqref{FirstBoundNablaSqrtuOn(0,T)} and \eqref{BoundOfEnergyFunctional} hold for  $L_1:=\f2{c_0}\max \big\{\mathfrak{F}(0)+\f{|\Om|}{e},c_2\big\}$ and $L_2:=\max \big\{c_2,\mathfrak{F}(0)\big\}$.  
\end{proof}
\end{lemma}
After these preparations, we can give the {\bf Proof of Theorem \ref{MainTheorem}}:  

Since the local-in-time classical solution $(u,w)$ of problem \eqref{AuxiliaryMainsSystem} is such that $w$ also solves $w_t \leq \Delta w +f(u)$ in $\Omega \times (0,T_{max}),$  by using a representation formula and the third assumption in \eqref{assumptions} we get 
 \begin{equation*}
 w(\cdot,t)\leq e^{t \Delta}w_0+\int_0^t e^{(t-s)\Delta} u^\beta(\cdot,s) ds \quad \text{ in } \Om \text{ and  for any } t\in(0,\Tmax).
 \end{equation*}
As a result, \eqref{LpLqEstimateGradient0}  with $\varphi=u^\beta$  
and the Young inequality lead to
 \begin{equation}\label{LInftyBoundW}
\begin{split}
&\lVert w (\cdot,t)\rVert_{L^\infty(\Omega)} \leq \lVert e^{t \Delta}w_0 \rVert_{L^\infty(\Omega)}+\int_0^t  \lVert e^{(t-s)\Delta} u^{β} \rVert_{L^\infty(\Omega)}ds  \\ & 
\leq  \lVert w_0 \rVert_{L^\infty(\Omega)}+C_S\int_0^t (1+(t-s)^{-\frac{\beta}{2}})\bigg(\int_\Omega u^2\bigg)^\frac{\beta}{2}ds \\ &
 \leq  \lVert w_0 \rVert_{L^\infty(\Omega)}+\frac{C_S}{2}\int_0^t(1+(t-s)^{-\frac{\beta}{2}})^2 ds
  +\frac{C_S}{2}\int_0^t \bigg(\int_\Omega u^2\bigg)^\beta ds \\ &
 \leq  \lVert w_0 \rVert_{L^\infty(\Omega)}+\frac{C_S}{2}t+\frac{C_S}{2(1-\beta)}t^{1-\beta}+\frac{2C_S}{2-\beta}t^\frac{2-\beta}{2}\\ &
 \quad 
 +\frac{C_S}{2}\beta \int_0^t\int_\Omega u^2+\frac{C_S}{2}(1-\beta)t,\quad t\in(0,T_{max}).
 \end{split}
 \end{equation}
From the hypothesis $\lVert v_0 \rVert_{L^\infty(\Omega)}<\chi^\frac{1}{\gamma-1}$, Lemma \ref{LemmaControllinguLoguAndnablauOveru} implies relations \eqref{FirstBoundNablaSqrtuOn(0,T)} and \eqref{BoundOfEnergyFunctional}; in particular,  \eqref{FirstBoundNablaSqrtuOn(0,T)} and  \eqref{Bound_of_u} make that  inequality \eqref{GagliardoNirenberg_u^2_on(0,T)} of Lemma \ref{lem:estimate:L2:spacetime} can be applied with the choice $\varphi=u$: hence
\eqref{LInftyBoundW}  
becomes
 \begin{equation*} 
 \begin{split}
 &\lVert w (\cdot,t)\rVert_{L^\infty(\Omega)} \leq  \lVert w_0 \rVert_{L^\infty(\Omega)}+\frac{C_S}{2}t+\frac{C_S}{2(1-\beta)}t^{1-\beta}\\ &
  \quad +\frac{2C_S}{2-\beta}t^\frac{2-\beta}{2}
  +\frac{C_S}{2}\beta \mathcal{K}(1+t)+\frac{C_S}{2}(1-\beta)t
  \leq L_{3}(1+t),\quad  t<T_{max},
  \end{split}
  \end{equation*}
where $
L_{3}:=\f{C_S}{2(1-\beta)}+\f{C_S\mathcal{K}\beta}2+\f{2C_S}{2-\beta}+\max\{\norm[L^\infty(\Om)]{w_0},\f{C_S}2+\f{C_S}2(1-\beta)\}.$  We now conclude the proof by invoking Lemma \ref{LemmaToEnsureTisInfty} with $C=\max\{L_1,L_2,L_3\}$.
 \qed 
 \subsection{Miscellaneous}  Even though for $β=1$ Young's inequality is not applicable in \eqref{GagliardoNirenberBoundu^beta+1}, for $c_0$ defined in Lemma \ref{LemmaControllinguLoguAndnablauOveru} relation \eqref{InequalityForF} is transformed  in 
\begin{equation*}
\frac{d}{d t}\mathfrak{F} +(c_0-C_{GN}^4m)\int_\Omega \frac{\lvert \nabla u\rvert^2}{u}\leq4C_{GN}^4m^2 \quad \textrm{on}\quad  (0,T_{max}),
\end{equation*} 
so that the remaining proof of Lemma \ref{LemmaControllinguLoguAndnablauOveru} continues to be valid if $m<\frac{c_0}{C_{GN}^4}$. Hence global existence of small-mass solutions to \eqref{problem} are recovered even for $\beta=1$: 
\begin{corollary}\label{MainTheoremBoundedness}
Let \eqref{assumptions} with $\beta \equiv 1$.  Then for any $\chi>0$ it is possible to find a constant $C_{GN}=C_{GN}(\Omega)>0$ such that for any $(u_0,v_0)$ as in \eqref{initdata} satisfying 
$\lVert v_0 \rVert_{L^\infty(\Omega)} < \chi^\frac{1}{\gamma-1}$ and $m=\int_\Omega u_0<\f{1-\chi \lVert v_0\rVert^{1-\gamma}}{2C_{GN}^4}$, there exists a unique pair of functions $(u,v)$ as in \eqref{solutionregularity}  which solve problem \eqref{problem}.
\end{corollary}

\end{document}